# Werner DePauli-Schimanovich
# Institute for Information Science, Dept. DB&AI, TU-Vienna, Austria
## Werner.DePauli@gmail.com

# Naïve Axiomatic Class Theory: A Solution for the Antinomies of Naïve Mengenlehre, K50-Set6

## 1 Abstract

Since the axioms in (Consi-CoS) are not recursively enumerable, NACT* is no axiom system in the classical sense[1]. Therefore we construct a series of partial systems which form a recursive axiom system too. Starting with the "dichotomic" systems NACT# and its variant NACT#4, we are going on to the "disjunctive" systems NACT+ and NACT+4, and eventually to NACT+Strat. After that we discuss the medium classes of these systems. Finally we present the inconsistent NSA-systems based on Not-SelfApplicability and explain their help for computational set theory.

## 2 Introduction

The system NACT* (= Naïve Axiomatic Class Theory Star) has been discussed in the papers "A Brief History of Future Set Theory" (K50-Set4) and "The Notion 'Pathology' in Set Theory" (K50-Set5b). It can shortly be described in the following way: Every formula of the collection of all hereditary-non-pathological formulas (with the only free variable "x") substituted into the class-operator generates a set.

NACT* is mainly interesting from a philosophical point of view and not from a logical one. Logicians have their ZFC and that will be sufficient for them. For mathematicians NACT* is not so important, because they still usually work with Naïve Mengenlehre despite the existence of ZFC. Nevertheless it might be interesting for mathematical logicians to construct a series of systems "converging" to NACT*. The first one are the dichotomic systems of NACT#. They are based on the first Patho-Thesis (1JT) that not both, "A" and "non-A", can be pathological[2]. This axiom system is very weak and without additional axioms it can even not be proved that $V \neq \{0\}$, which is a model of axiom (5) and (6). Therefore we must add 4 axioms which are similar to ZF-axioms to form NACT#4.

The NACT#-Systems produce only small sets (and the complements of the proper-classes). This does not correspond to our intention to allow also universal sets. Therefore we have to throw away the dichotomy-law and replace it by our normal disjunction-law (which always includes the case "or both"). With this new disjunction-law we construct the system NACT+ which allows also universal sets. But also here, the pure frame "NACT+" is too weak. We cannot even prove V ≠ Omega, which is a model of (5a) and (6c). Therefore we have to add the 4 ZF-like axioms again to get NACT+4. That is the easiest way. Finally NACT+Strat is defined, a combination of NACT and NF.

---

[1] For (Consi-CoS) and NACT* see "A Brief History of Future Set Theory" and "The Notion 'Pathology' in Set Theory", both articles in this book.
[2] See "The Notion 'Pathology' in Set Theory", paragraph 6.

After setting up the axiom systems, some properties of sets of these systems are discussed. E.g. the "slim", "mightly" and "medium" classes are investigated and the question is considered, under what circumstances the medium classes are sets.

## 3 Axiom-Systems of NACT starting with NACT# and NACT#4

The goal of the construction of NACT# and its extensions is to find an axiom system which can (at least partially) blow up its own cumulative hierarchy of sets. And to generate as much sets as possible. Class-operator, Church-schema, Extensionality and (AC) together we want call again the frame of the class theory CT (usually NBG). In addition to the 4 axioms of the CT-frame, NACT# uses the two special axioms (5) and (6) which yields altogether:

(1) class-operator {|x: A(x)|} := {|X: Set(X) & A(X)|},
(2) the Church-Schema (CS) as comprehension scheme for classes,
(3) extensionality axiom as usual,
(4) a suitable axiom of choice (it may be (DC) or the super-strong Ord ~ UC),
(5) instead of the Zermelo-Fraenkel axioms as the only axiom, the "general set generator":
(Set-nonEquiSetKo) := [Set(X) <=/=> Set(ko(X))],
(6) the "ordering" axiom:
(SlimSet) := [Slim(X) ==> Set(X)], where Slim(X) := [card(X) < card(Ko(X))].

Because of (5) we call the NACT#-systems also "dichotomic" systems.
NACT# is too weak to generate enough sets. Therefore we have to add 4 axioms, thus yielding NACT#4.
(#1) Slim(Omega),
(#2) Slim(X) ==> Slim (Power(X)),
(#3) Slim(X) & [forall y in X: Slim(y)] ==> Slim(UX), where "U" is the large union.
(#4) Slim(X) & Function(F) ==> Slim(F[X]), where F[X] is the image of X under F.

Restricting NACT# or NACT#4 to Mirimanoff sets (without an infinite element-sequence) -- in the same way we restricted the classes of NBG to sets -- we can construct a well-founded system NACT° which could be acceptable for the ZF-people.
    Establishing our next system we are making a step forward towards NACT*.

## 4 The "disjunctive" axiom systems NACT+, NACT+4 and NACT+Strat

If the reader has already accepted that NACT#4 (and NACT° respectively) are (empirically) consistent (e.g. relatively to NBG) then we can proceed to the next system NACT+.

Axioms (1) to (3) (i.e. class operator, comprehension scheme and extensionality) are the same as in NACT#, but for (4) we prefer a weak axiom of choice. In (5) from NACT# we replace the non-equivalence by a disjunction, and thereby getting

(5a) := (SetDisKoSet): = [Set(X) or Set(ko(X))].

As ordering axiom (6) we use instead of (SlimSet) of NACT# the axiom (SlimMightySet) which generates the "mighty" sets or "anti-sets" too (in addition to the slim sets). We use this new axiom in the (equivalent) formulation
(6c) := (SlimSet&KoSet): = [card(X) < card(ko(X)) ==> Set(X) & Set(ko(X))] .

If we add to NACT+ the 4 axioms (#1) to (#4) we get NACT+4.

   For NACT+ it would also be helpful if we add Quine's Rule of Stratification of the system New Foundations forming in such a way NACT+Strat. The stratification rule (StratSet) would be formulated here as axiom (StratCoS):
Stratified(A) ==> Set({|x: A(x)|}).

Because in NACT+Strat already all slim classes are sets, here the problems of NF with ordinals, omega, etc, do not arise. By the addition of stratification the 4 axioms (#1) to (#4) become redundant. But it may simplify the new system NACT4+Strat, if we still use them.

   Until now no investigations had been made to compare NACT#4, NACT+4 and NACT+Strat, but it would be certainly an interesting endeavour.

Now let us go on with the investigation of our question about ordering axioms for the medium classes.

## 5 The medium classes of the NACT axiom-systems

After showing his famous proof of the incompleteness of arithmetics, Kurt Gödel said: "And now Mengenlehre". 1938 he proved with the method of inner models (of constructible sets) the relative consistency of (AC) and (GCH) with ZF. 1963 Paul Cohen showed with the forcing method the corresponding fact, that also ⌐ (AC) and ⌐ (GCH) are relative consistent to ZF. These proofs had been very difficult already for ZF. But similar proofs for NACT#4, NACT+4 and NACT+Strat will be certainly much more difficult. The medium classes make some difficulties because of their "paradoxical" behaviour.

We have first to define:
Slim (X): = card (X) < card (Ko(x)) and
Mighty (X) : = card (X) > card (Ko(x)),
Slim : = {|x : Slim (x) |} and Mighty : = {|x : Mighty (x) |},
Medium : = {|x: Medium(x)|}, where Medium(X) : <==> X ~ Ko(X).

We see immediately that Slim u Mighty = Ko(Medium). First for the moment we assume that set(Slim) is valid and also set(Mighty), and therefore proper-class(Medium). We do not like very much for Slim(Slim) to be valid and therefore also not "Slim in Slim". (But since already in this case "{Slim} in Slim", NACT# and NACT+ are in any case incompatible with the foundation axiom.) It would be more in line with our intuition if non-Slim(Slim) or (what is the same) Medium(Slim) and therefore Slim ~ Ko(Slim). Fortunately Slim(Slim) implies

Cantor's Antinomy, and therefore Medium(Slim) is valid. A similar argument can be applied to the class Mighty, but Medium is different.

For x in Slim u Mighty (where "u" is the union) the complement exists and is again in Slim u Mighty. Similarly we want consider MediumC : = {|x: x in Medium & ko(x) in Medium|}. Since for the complements of the pathos the complement does not exist, I want to establish the following conjecture:

„Patho = MediumNC", where Medium NC : = Medium \ MediumC.

If this conjecture can be proved it would be a first set-theoretical characterization of Patho. A partial characterization can be the following formula:

non set (X) ==> X ~ Ko(X) & set (Ko(X)) & Ko(X) non-in Ko(X).

"Proper-classes are complements-equipollents and complements-sets and complements-not-themselves-containing."

or the contraposition (NACT*):

X ∤ Ko(X) or non set (Ko(X)) or Ko(X) in Ko(X)  ===>  set(X).

"Sets are these classes which are not-co-equipollent or co-properclasses or co-selfcontaining."

With other words: "The co-not-equipollents or the co-properclasses or the co-selfcontainers are sets!"
I call this the (NACT*-conjecture).

For NACT+Strat we can see that (StratCoS) mainly influences the set-generator for the class Medium because, with A, also non-A is stratified, and therefore both are sets now. Hence NACT+Strat contains the slim sets, the mighty ones, the stratified mediums and the complements of the pathological classes. That's more than sufficient.

Because axiom (5a) leaves open whether one or two classes of the unstratified part of Medium are sets, NACT+Strat does not generate a unique maximal universe like NACT° or at least a "nearly unique" maximal universe like NACT# or NACT#4. It would be fine if we could add axioms to NACT+Strat to force it to generate a (nearly) maximal universe in an absolute sense, i.e. that nearly every class (except the pathos and sub-pathos) that can be assumed to be a set (consistent with the already created ones) is already in the universe.

Since in NACT# of every pair of formulas A and non-A, exactly one creates a set, NACT#4 generates "Maximal Possible Universes" of "small" sets, which can consistently coexist together. The small sets are here the slim classes, and classes which are provably the complements of proper classes. But since it can also happen that "X ~ Ko(X)" is valid for a set, some  "medium" classes (with X ~ Ko(X)) can also be small (in this special sense of NACT#4). Unfortunately we do not know which one of the X and Ko(X) is the set (and both cases maybe possible). Therefore the maximal (consistent) universes are not unique.

The set-generator axiom (5) of NACT# only guarantees that either a class or its complement is a set. But it does not say which one of the pair of classes is the set and the ordering axiom (6) regulates this only for the slim classes. Therefore the mighty classes with "card(X) > card(Ko(X))" are proper classes in NACT#. But what happens with the whole collection of medium classes "X ~ Ko(X)"? Are we able to find an ordering axiom for all the medium classes? Or a list of axioms for parts of the medium classes? If we restrict the universe of NACT#4 to non-medium sets, then it is unique.

This is the point where we started to construct the systems NACT$ and NACT*[3]. Important for logicians is the fact that NACT* does not fall from heaven but rather is the limit of a sequence of systems of NACT. Since NACT#4, NACT+Strat and NACT* imply the (parameter-free) axioms of NBG (whereby some of them are restricted to a slim, well-founded, Cantorian, or Miri domain) and since NACT+Strat and NACT* imply the parameter-free NF (probably with the Axiom of Counting), NACT* is the strongest known system of all, generating a maximal universe of sets in an absolute sense.

## 6 The inconsistent systems NACT-PriN(G)SA(2) and NACT-SiN(G)SA(2)

During my work on the Naïve Axiomatic Class Theory NACT I investigated also a series of 8 elegant systems I called short NSA (or better NACT-NSA). PriNSA is the weakest one of these 8 of similar set-theoretical systems based on the (normal) "Not-SelfApplicability" NSA.

Let CT be the class-theoretical frame as usual consisting of class-operator, Church schema, extensionality and (AC) with well-founded parameters.

Call a wff A self-applicable or SA, if A is valid for the class generated by it:
SA(A) := A({|X: set(X) and A(X) |}).

Then the only set-generator of system (PriNSA) is:
not SA (A) $\implies$ set ( {| x: A(x) |}).

All together there are 4 similar systems based on the NSA-property:

PriNSA "Principle of Not-allowed SelfApplication"
not A({| X: set(X) and A(X) |}) $\implies$ set ({| x: A(x) |}).

SiNSA "Strong rule of NSA"
not A({| X: set(X) and A(X) |}) $\iff$ set ({| x: A(x) |}), where A does not contain "set(X)".

PriNSA2 "2-fold Principle of NSA"
not A({| X: set(X) and A(X) |}) $\implies$ set ({| x: A(x) |}) & set ({| x: non-A(x) |}).

and SiNSA2 "2-fold Strong rule of NSA"
not A({| X: set(X) and A(X) |}) $\iff$ set ({| x: A(x) |}) & set ({| x: non-A(x) |}),
again with the condition, that A does not contain "set(X)".

---

[3] See „A Brief History of Future Set Theory" and „The Notion ‚Pathology' in Set Theory".

If you define a Generalized SelfApplication GSA and its negation NGSA, you get 4 other systems: PriNGSA, SiNGSA, PriNGSA2 and SiNGSA2.

The definition of GSA can look like the following way:

A*[ {|x: A(x) |} ] := A({|x: A(x) |}) or [exist n>=1: exist x1, ... xn: {|x: A(x) |} in x1 & x1 in x2 & ... & A(xn)].
That means more or less that the class generated by A is in the element-closure of A.
     SiNGSA and SiNGSA2 can only be used together with the meta-rule (Meta-SiNGSA(2)):
If SiNGSA(2) entails proper-class({|x: A(x) |}), or if we assume this, {|x: A(x) |} cannot be substituted into SiNGSA or SiNGSA2.
     It would also be good to apply this meta-rule to PriNGSA(2) too, because it supports our intended meaning of GSA. [i.e.: since proper-class({|x: A(x)|}) makes the $2^{nd}$ term in GSA always false, it reduces GSA to SA.]
     The new systems are completely different than the old ones. e.g. "Si" := "the class of all singletons" is a set in SiNSA and contains its own singleton, but it is a properclass of SiNGSA.

It is very interesting how the machinery of these 8 systems works and to compare them. This is the reason why I presented the systems here to the reader. Unfortunately it turned out that even the weakest system is inconsistent. But you can compute the degree of inconsistency, what is still an appropriate question for computational logicians.[4]
     When we work with PriNSA, essentially 4 different cases can happen:
(1) PriNSA(Ai) is inconsistent, i.e. it generated a contradiction together with other substitution-instances of Ai. In this case no proposition can be made about the class generated by Ai [i.e. if it is a set or a proper-class] which cannot be made already by the CT-frame [e.g. that Ru is a proper-class].
(2) NSA(Ai) in PriNSA is not decidable. Then again no proposition can be made about class Ai [which cannot be made already by the CT-frame].
(3) SA(Ai) is valid. Again no proposition can be made by PriNSA. And
(4) NSA(Ai) is valid. Then class Ai is a set.

## 7 Properties of the system PriNSA
PriNSA does not work like a Hilbert-calculus but like a production system PS producing wffs Ai with increasing length. It is more practical if we assume that the wffs Ai are parameter-free, i.e. with the only free variable "x". PS produces also the terms t with increasing length [with only bounded variables]. Therefore we can check for one wff Ai after the other if PriNSA(Ai) is inconsistent with respect to the already done derivations. If this turns out for a special Ai, PriNSA(Ai) is not allowed going on to use Ai and any wff Bj which contains Ai. All further derivations from Ai (and the already derived facts from it) are forbidden. This "Wittgenstein-Rule"[5] is of course a meta-rule outside the deduction calculus of logic. We will describe this procedure later in more detail.

---

[4] See Schimanovich [1981].
[5] See Wittgenstein [1956].

Proceeding this way we can also compute the consistency-degree of an axiom-system.[6] For computational logicians it has also sense to investigate inconsistent axiom-systems like Arrow's paradox[7], the Naïve Mengenlehre, or the NSA-systems. E.g. if the limes of fraction [between the number of derived false formulas and the number of all wffs] is 0 [using the Wittgenstein-stoprule], then we call this system nearly-consistent. This way we can also formulate a partial-consistency.[8]

Because of this interest of the computational logicians on the NSA-system, let me here list up some facts which are already derivable after some steps, long before that the first contradiction can be derived.

For comparison we want to note the special situation of the complement Ko(Ru) of the Russell-class Ru as an example.

Case 1: In PriNSA both is possible: "ko(Ru) in ko(Ru)" and also "Ko(Ru) non-in Ko(Ru)". In all 4 NSA systems ko(Ru) is a set, and ko(Ru) and its derivatives are the only complements of pathological classes, which are sets in these systems.

Case 2: In SiNSA ko(Ru) non-in ko(Ru) is valid.

Case 3: In PriNSA2 ko(Ru) in ko(Ru) is valid.

Case 4: In SiNSA2 ko(Ru) in ko(Ru) is also valid.

This formulation of SelfApplication allows us to establish a system NACT-PriNSA where neither the universal-classes, nor the Pathos or its complements are sets except ko(Ru). (Therefore this system is very similar to ZF and seems to me very important.)

The antecedence of PriNSA is the (normal) negation of the "SelfApplicability" SA (or s.a.). It is a formula-property like the stratification (but not so easy as this) or like Hereditary-non-Patho (but not so difficult as this). Therefore we want to define SA in the following way: SA(Ai) is the SelfApplicability of Ai({|X: set(X)&Ai(X)|}) and it says that the class {|x: Ai(x)|} fulfills the property of Ai. Not Ai({|x: Ai(x)|}) says of course that it should not fulfill it: NSA(Ai).

E.g.: $\emptyset = \{|x: x \neq x|\}$ is not selfapplicable, because $\emptyset \neq \emptyset$ is wrong, but $V = \{|x : x = x|\}$ is selfapplicable, because $V = V$ is true. Also Ru = {|x: x non-in x|} and its derivatives are SA, because Ru non-in Ru, etc. The same is valid for On, Cn, Fund, Miri, etc. But also the complement Ko(On) of On = {|x: non On(x) |} and its derivatives are SA, because it is no ordinal-class, etc., and also the complements of the other Pathos and their derivatives too. Only Ko(Ru) and its derivatives are exceptions, because it is not clear if ko(Ru) in ko(Ru) [and therefore SA] or the contrary, etc. If ko(Ru) non-in ko(Ru), it is not SA and therefore a set [and the same is valid for the derivatives of ko(Ru)].

---

PriNSA [or better NACT-PriNSA, because PriNSA is embedded in the CT-frame] is the weakest system of these 4, but it still contains also not-wellfounded sets. Therefore we should restrict the sets to wellfounded ones [like we restricted the classes to sets], i.e.:

{x: A(x)} : = {| X: Fund (X) & set (X) & A(X) |}.

Note that here the first class-operator has no vertical bars after and before the set-brackets/parentheses. We want to call this system of restricted classes NACT-PriNSA$^0$, which is similar to ZFC$^0$.

## 8 Properties of the system SiNSA

If we strengthen the implication in PriNSA and replace it by an equivalence, it yields the "Strong rule of Not-allowed SelfApplication" (SiNSA):

Not A({| X: set(X)&A(X)|}) < = = > set({| X: set(X) & A(X)|}).

A should not contain the set-predicate set (…X…) with any function of X in the NSA-condition of SiNSA. Only implicit in bounded set-variables, i.e. in "exist X: set(X)&B(X)" and "forall X: set(X) => B(X)". If we expand free set-variables, the antecedence "set(X)" will stand before the formula A. And about the case that there are other variables Y, Z, etc. [instead of X] we cannot make a final decision until now.

   Thus our system NACT-SiNSA consists of (CT1) to (CT4) and (SiNSA).

NACT-SiNSA implies the following ZF-axioms:

(ZF1)  set (Omega)

(ZF2)  set (Power (b)), where b is a free parameter,

(ZF3)  Fund (b) ==> set (Ub), where "U" is the large union, and

(ZF4)  Fund (d) ==> set (F[d]), where F[d] is the class of all function-values F(x) for x in domain d.

(ZF5)  set ({|a, b|}), where a and b are free parameters.

Since Omega = {|x: natural number x|} is no natural number, it is not SA. Neither the pair-class {|a, b|}, nor the power-class of the set b are SA, because the pair cannot be a or b, and the power cannot be a partial-set of b. [This is trivial for the well-founded parameters, but probably also valid for not well-founded sets.] If b is well-founded it is also valid for Ub, which is therefore a set. If the domain d is well-founded, the image F[d] is no function-value of F and therefore not SA and a set. So we can see that the results of all 5 ZF-axioms are sets because of SiNSA.

SiNSA probably implies also Tarski's axiom of exorbitant sets (with well-founded parameters) or the strong cardinality axioms (like the axioms of measurable cardinals, super-compact cardinals, Woodin cardinals, etc, also with well-founded parameters) and is therefore much stronger than ZF. SiNSA is recursively enumerable but probably not decidable. We do not care about that, because SiNSA is a unique principle to generate sets in contrary to the arbitrary chosen axioms of ZF.

   But we should also confess here that SiNSA produces some strange sets too. [Strange because we are educated to think in ZF-termini. But these strange sets are in no case paradox and not at all pathological!] E.g. the mentioned ko(Ru) and its derivates are strange sets. But also the set of all singletons si := {| x: exist y: x={|y|}|} and its derivates. si contains not only {|ko(Ru)|} but also {|si|}.

Therefore also not-wellfounded sets exist in NACT-SiNSA. Because of these strange sets it can be the case that it will turn out that also (ZF2) and (ZF5) need the prefix "Fund (a,b) ==>" too. [Until now the author did not find an example which entails this.] Also for (ZF4) it may turn out that we have to add "Fund(F[d])" [or replace "Fund(d)" by it]. In ZF our (ZF5) follows from unrestricted (ZF3). But since (ZF3) needs a restriction [because Usi would produce V as a set], while (ZF5) not, we listed it separate.

Restricted (ZF4) implies only restricted separation. This is good, because unrestricted separation would generate (si intersection Ru2) [where Ru2 := {|x: x non $in^2$ x|}] as a set, what turns out to be a sub-class of Ru2, which is probably a proper-class too. Note that in NACT-SiNSA no selfcontained sets x [with x in x] exist, and ko(Ru) is therefore = 0. If (AC) would not be restricted to wellfounded parameters, the choice-function with domain si would select the universal proper-class V as choice-set v.

Since in $ZFC^0$ the conditions "Fund (…) ==>" can be dropped, NACT-SiNSA implies $ZFC^0$, in the sense that all sets which can be produced in $ZFC^0$, exist also in NACT-SiNSA. Therefore:
$ZFC^0$ = < NACT-SiNSA := forall x: [setX in $ZFC^0$ ==> set X in NACT-SiNSA],
V sub $ZFC^0$ = < V sub NACT-SiNSA and probably also
V sub ZFC = < V sub NACT-SiNSA are valid.

## 9 Properties of the systems PriNSA2 and SiNSA2
NACT-SiNSA may be interesting for the ZF-people, but it does not allow universal sets, which was our original goal to establish systems that allow them. Therefore we consider PriNSA2:
not A({|x:A(x)|}) ==> set ({|x:A(x)|}) & set ({|x: non A(x)|}).
    PriNSA2 is similar to PriNSA3:
not A({|x:A(x)|}) or not non-A ({|x: non-A(x)|}) ==> set ({|x:A(x)|}).
    With these 2 principles we get NACT-PriNSA2 and NACT-PriNSA3. But we are not very interested in NACT-PriNSA3 and leave it to the reader to check NACT-PriNSA2.

Furtheron we want again strengthen the implication in PriNSA2 to an equivalence. By this change we get SiNSA2:
not A ({|x:A(x)|}) <==> set ({|x:A(x)|} & set ({|x: non A(x)|}).
Here A again should not contain the set-predicate set (…X…) of some function of X.
    This 2-fold strong rule together with the CT-frame forms the system NACT-SiNSA2.
It is also possible to consider NACT-SiNSA3 [which has nothing to do with NACT-PriNSA3]. It uses SiNSA3:
not A ({|x:A(x)|}) ==> set ({|x:A(x)|}) & set ({|x: non A(x)|}), and
A ({|x:A(x)|}) ==> non set ({|x:A(x)|}) or Mighty ({|x:A(x)|}).
Mighty (x) means card (x) > card (Ko(x)).

But SiNSA2 is better and has more symmetry-properties. Eg.: with 0, {0}, Omega, etc., also V, V\{0}, V\Omega, etc. becomes a set, also {v}, v\{v}, etc.

But if a class is SA like the universal sets or like the Pathos, it or its complement becomes a proper-class, or in SiNSA3 Mighty [except ko(Ru)]. So NACT-SiNSA2 consists of the CT-frame and SiNSA2.

But fortunately it does not imply all unrestricted ZF-axioms (ZF1) to (ZF5), but only the ones restricted to small sets. As representative for the property "small" we will use this time the terminus "Fund" and define:

(Fund)      Fund(x) : = exist y in X : y intersection X == empty.

With this definition we formulate:

(ZF1R) = (ZF1) set (Omega),

(ZF2R) Fund (b) ==> set (Power(b)),

(ZF3R) Fund (b) ==> set (Ub), and

(ZF4R) Fund (F[d]) & Function (F) ==> set (F[d]).

This time we drop the pair-set axiom, because it follows from (ZF3R).

Because of the restriction of (ZF2R), e.g. Cantor's antinomy is not derivable. But also if we allow Power(v) to be set, the "across-element" q = {|x: x non-in F(x)|} for all F : v $\rightarrow$ Power(v) need not to be set.

Since the identity I is a function v $\rightarrow$ Power(v), q degenerates to Ru. Because of the restriction of replacement, which yields a restricted separation, set(Ru) is not derivable, etc. Since with b also Power(b) is Fund, it is not SA [i.e. no subclass of the well-founded b] and therefore a set, and also Ko(Power(b)). A similar reasoning shows that NACT-SiNSA2 also implies (ZF3R) and (ZF4R), and therefore produces all sets, which ZFC$^0$ can produce.

NACT-SiNSA2 is similar to the Syntho-Mengenlehre in Schimanovich [1971a] and K50-Set1 or in Schimanovich [1971b] and K50-Set2, which is ZFC with restricted replacement and (AC) [and an Omega-Existence-Axiom instead of the axiom of infinity], plus the Complement-Existence-Axiom, and therefore more important than ZFC.

It seems that the definition of SA we gave before is sufficient. But later it can happen that we have to restrict the wffs A to parameter-free formulas. This helps us to show that ZF-axioms are valid. Also it may become interesting to generalize SA by defining the generalized selfapplication GSA [i.e. the SA in the wider sense]

A*[{|x: A(x) |}] : = A({|x: A(x) |}) or [exist n ≥ 1: exist $x_1$, $x_2$, …. $x_n$: {|x: A(x)|} in $x_1$ & $x_1$ in $x_2$ & … & A($x_n$)].

That means more or less that {|x: A(x)|} is a member of the element-closure {|x: A(x)|}* and the last link in the chain of containing elements has the property A. To remember on this fact we write the class {|x: A(x)|} in corner-brackets "[" and "]" and put a star in front of it after the A. In future we will call this generalized selfapplication GSA. But the systems with PriNGSA(2) or SiNGSA(2) are different from these without G. [E.g.: Si is a proper-class in these systems.]

Since we know that NACT-PriNSA, NACT-PriNSA2, NACT-PriNGSA, NACT-PriNGSA2, NACT-SiNSA, NACT-SiNSA2, NACT-SiNGSA or NACT-SiNGSA2 are inconsistent, we want to concentrate ourself onto the still interesting question, after how many steps (concerning to a standard algorithm of the production system PS producing wffs) "falsum" will be derivable. One can also define the degree of inconsistency as limes of the fraction (number of false wffs) through (number of wffs) produced by PS. That is still an interesting problem for computational logicians. (See Schimanovich [1981]).

## 10 NACT# and NACT+ combined with PriNSA or SiNSA

We want now to try another way to enlarge NACT# by adding the mentioned "Principle of Not-allowed SelfApplication", but this time for parameter-free wffs A, which is a weakened form of PriNSA:

$PriNSA^0 := $ non $A^0(\{| \; x: A^0(x) \;|\}) ==> $ set $(\{| \; x: A^0(x) \;|\})$.

$A^0$ should mean that A is nearly-closed, i.e. has no other free variables than x.

If (under the assumption that $\{|x: A^0(x) \;|\}$ is a class the parameter-free formula $A^0$ does not apply to itself, then class $\{|x: A^0(x) \;|\}$ must be a set. (In future we will drop the parentheses "(" and ")" as superscript, and write only A. Onto the Pathos like the class of all ordinals or cardinals, etc. it has no influence, because the Pathos are self-applicable. They can not become a set and generate a contradiction. But $PriNSA^0$ also does not support the construction of universal sets. That is good because in NACT# universal classes are no sets. But for our considerations here we can use it to form the system $NACT\#PriNSA^0$ consisting of the CT-frame (1) to (4), the set generator (5) and (6b) : = $PriNSA^0$ for parameter-free wffs.

NACT#PriNSA : = (1) & … (5) & $PriNSA^0$ fulfills (#1) and for all closed terms X also (#2) – (#4). By meta-mathematical reasoning over the number of axioms [especially PriNSA], we can show that (#2) – (#4) are valid for slim stets X. But e.g. (#2) or (#4) are not valid for X = V, (also if V would be a set v). Therefore this system blows up its own cumulative hierarchy and therefore the axioms (#1) to (#4) are not longer necessary. And (6) is also not needed [because PriNSA implies it]. And because of (5) the universal-like classes cannot form a set.

If necessary we should again generalize the SelfApplicability SA. This new generalized GSA A* [{| x: A(x) |}] can be formulated including the case of forbidden infinitely-long descending element-sequences too, what yields MSA the corresponding principles are PriNMSA and SiNMSA. But we want note that already PriNSA does not generate the complements of Miri, On, etc. as sets [because they are as well selfapplicable as the Pathos]. This is done only by (5). We mentioned already that PriNSA together with NACT# forms NACT#PriNSA. All classes with parameter-free formulas $A^0$, for which $A^0(\{| \; x: A^0(x) \;|\})$ is valid, are sets and by dichotomy law (5) therefore their complements are proper-classes.

E.g.: since "x ≠ x" is not A($\{| \; x: A(x) \;|\}$), 0 is a set and V a proper-class. The same is valid for all constructions "from the bottom" and their complements. Because the Russell-class Ru and its derivatives are proper-classes (and of cause pathos), ko(Ru) and the complements of the derivatives are sets [concerning to (5)]. In fact you cannot decide (by logical means only) if ko(Ru) in ko(Ru) or not. Since ko(Ru) in ko(Ru) would be a selfapplication, ko(Ru) non-in ko(Ru) is better to be valid, and the same for the derivatives. But the complements of the other Pathos are selfapplicable [and therefore PriNSA is not applicable]. So the reader can see that NACT#PriNSA forms already a very strong system.

Similar to NACT-PriNSA2, we want now construct NACT+PriNSA2 by adding the parameter-free PriNSA$^0$2 to NACT+. To remind the reader, we repeat this "2-fold Principle of Not-allowed SA for parameter-free wffs":

not A$^0$({| x: A$^0$ (x) |}) ==> set({| x: A$^0$ (x) |}) & set({| x: non A$^0$ (x) |}).

In future we will drop again the superscript "$^{()}$" after A.

NACT+PriNSA2: = (1) & … (4) & (5a) & (6c) & parameter-free PriNSA2.

If we use GSA (or MSA) we get the corresponding systems.

We can now discuss NACT+SiNSA2, which is much stronger then SINSA. [More than that: It generates also the complements of small sets as sets.] Therefore we could not use it for NACT# any longer and moved for this reason to system NACT+. [Here we use the set-generator (5a) and ordering axiom (6c).] As before we have

NACT+PriNSA3:= (1) & … (4) & (5a) & (6c) & parameter-free PriNSA3.

But this system does not bring more than NACT+PriNSA2, and is not so elegant. It is therefore better to investigate NACT+PriNSA2 or NACT-SiNSA2.

## 11 Not SelfApplicability as Additional Condition and Computational Set-Theory

A lot of facts we can derive in few steps from PriNSA(2) and SiNSA(2) and also from their combinations with NACT# and NACT+. But now we are already at the point where the first inconsistency appears.

Consider si: = {| x: exist y : x ={y}|} and {| si |} and you will find everything o.k. But modify this to si' : = {| x: exist y : x = {y} & {y} non-in y |}. Then:

{si'} in si' <==> {si'} non-in si'.

A clear contradiction![9]

The class-comprehensing formula here is a little bit longer than the other formulas we studied until now. Therefore it is the first inconsistency the production system PS generated. This implies that the consistency-degree of the NSA-systems is probably 1.

This Russell-like paradox can be derived because we did not demand the "Stratification" property in addition to NSA. It would be very interesting to investigate the following schema Strat&NSA:

Strat (Ai) & not SA (Ai) ==> set ({| x: Ai(x) |}), where Ai does not contain "set (…X…)".

From Ronald Björn Jensen's famous article we know that NF is relative consistent to ZFU. The addition of NSA forces us to produce less sets than in NF. To compensate this and to fulfill the symmetry of the stratification concerning the complement we want to use better Strat-NSA2:

Strat&NSA(Ai) ==> set ({| x: Ai (x) |}) & set ({| x: non Ai (x) |})

or short: StratNSA(Ai) ==> set ({| x: [non] Ai (x) |}), where set(X) $\not\subseteq$ Ai.

But the addition of NSA to stratification as precondition for the sethood does not make NF easier and does not bring us very much. Therefore it is better to choose another condition where NSA can be added. If we take Miri(x) [i.e. the property that there is no infinitely descending element chain in X] or Fund (X) [i.e. that X is well-founded], similar inconsistencies as before can be derived. [10] Therefore I suggest to choose Hereditary-non-Patho as second precondition for the sethood, what turns out to be more or less equivalent to the Wittgenstein-rule.

---

[9] I want to express my thanks to Randall Holmes who gave this counterexample to PriNSA to me.

But the combination of HnP and NSA is more efficient than HnP alone. Let us call this system NACT**, which is established by the combination HnP&NSA. A lot of proper-classes [and its derivatives] are already filtered out by NSA. But NSA alone is inconsistent and we need therefore also HnP. Remember: HnP(Ai) : =
Forall wff Bj partial-formula of Ai: from set ({|x: Bj(x)|}) is not-derivable falsum within the CT-frame.

This way we could save the inconsistent NSA-schemata from its insignificancy. But we can also use the combination HnP&NSA as a basis for a computational set theory. First we have to enumerate the wffs according to their "length": verum, falsum, x in x, negation, x in x1, x1 in x, x1 in x1, negation, quantification over x1, x1 in x2, x2 in x1, negation, quantification over x1 and x2, negation, quantification over x2 and x1, negation, conjunction of 2 preceding formulas, etc, forall x1: x in x1, forall x1: x1 in x, forall x1: x1 in x1, forall x1 forall x2: x1 in x2, forall x1 forall x2: x2 in x1, negation outside of the preceding formulas, … conjunctions of 2 of the preceding formulas, etc. One can easy find out the rules which generate this chain of formulas. Here we use the wffs only in expanded form, i.e. without explicitely writing down "exist", "or", "implies", "is equivalent" or "is equal". [We use (EE) as definition of equality.]

Once the enumeration of the wffs is established we can start the theorem prover. Of course we will stop it after 10 millions of derivation steps and probably we will work with time slices and swop. We shall try first if a formula Bj turns out to be inconsistent and then we drop all formulas Ai containing Bj. Usually this is easier to show than the contrary what may furthermore be undecidable. Therefore we should move to show NSA or SA of Bj. If after 10 millions of derivation steps we still do not have any result we consider HnP&NSA(Bj) as probably undecidable. Otherwise the class generated by Bj is a set. I want to note here that this procedure is not only a good basis for a computational set theory, but also the way in which mathematicians think, i.e. in terms of Naïve Mengenlehre!

## 12 Conclusion

In NACT all classes are predicate-extensions. (In NACT* the closed consistent predicate-extensions are even sets. And in the other systems of NACT only those classes are sets which satisfy the axioms.) In NAM[11], on the other hand, only "normal" sets are predicate-extensions. In NAM# e.g. abnormal sets are all identified with the universal set (shown by Cantor to be inconsistent for the old Naïve Mengenlehre).[12] With such technical tricks we can show that naïve set theory is not inconsistent after all.

---

[10] Miri and Fund are only equivalent to each other because of (AC), but the "equivalents of the (AC)" needs not to be longer equivalent in stronger set theories than ZF. Therefore we assume that in the systems we use here Miri and Fund are different.
[11] See "Naïve Axiomatic Mengenlehre for Experiments".
[12] See Halmos [1960].

Thus, in addition to the traditional Zermelo–Fraenkel set theory ZFC[13] with the axiom of choice[14] (or, alternatively, the von Neumann–Bernays–Gödel class theory NBG[15]) and Willard Van Orman Quine's New Foundations NF[16] (alternatively Mathematical Logic ML), we now have available NACT (and its modifications) as a third group of systems of set theory. And I am convinced that ZFC will gradually lose its dominant position — because ZFC "follows from" NACT, i.e. is "justified" by it. Set theory also serves as the foundation for all mathematics as well as concept formation in all empirical theories. Mathematicians will, I think, prefer NACT because it is nothing other than a formalization of naïve set theory, which they use routinely anyway. NACT justify this traditional use of naïve set theory.[17]

Nonetheless, Zermelo–Fraenkel still remains part of the current academic Zeitgeist — which was the "Ungeist der Zeit" for Kurt Gödel[18], the bad spirit of the times. It is for this reason that NBG as well is the dominant standard for foundational work in mathematical logic. But this role can be better played by NACT — I assume it will work properly (e.g. in investigations on measurable cardinals).

## 13 A philosophical Outlook
A quite important reason for research on NAM is that this represents a family of systems in which the notion of normal set is suitably characterized.[19] Mathematicians, who typically love to hover in their Platonic heaven, have a good sense for what normal sets are, and for that reason never blunder into inconsistencies when using naïve set theory intuitively. Inconsistencies arise only through incorrect symbol manipulation[20] with "pathological" sets.

ZFC, for example, defines normal sets only in part and implicitly through its axioms. There are, however, a large number of normal sets not comprised here. For that reason, researchers deliberate on whether to assume the existence of measurable cardinals or exorbitant ordinals or not, which clearly seem to be normal and unproblematic.

---

[13] See Suppes [1960] and Bernays-Fraenkel [1968] .

[14] See Rubin-Rubin [1978].

[15] See Gödel [1938].

[16] See Quine [1969].

[17] It does not matter that NAM* and NACT* are no axiom-systems in the classical sense because „Consistency" (or "Decidability") of wffs are no recursive enumerable properties. Since the ruin of Hilbert's Program we need no axiom-systems of set theory any longer where we can prove that large parts of math can be implemented in it. But we need still a justification for the contentual method of the mathematicians. That's the reason why I investigated these systems. And it is even imaginable that for theorem proving NACT works better than ZFC.

[18] See Schimanovich-Galidescu [2002] in the references.

[19] See DePauli-Schimanovich [2002a] und [2006a] in the references.

[20] Mathematicians and logicians think in models and the models are consistent (by definition). But inconsistencies result from a mistake in the derivation machinery. Therefore we have to study this machinery if we want to discover the occurrences of contradictions. (Only Frege and Russell thought in this way.) My point of view is that NACT and NAM can be interpreted as suitable explications for logicism.

Also NF determines the notion of normal set inductively by the comprehension scheme with stratified predicates in a well-known way. But Quine's stratification condition seems arbitrary, so it is not clear whether it really characterizes normal sets properly. In any case, it is not intuitive for most mathematicians.

Other set theoretical systems, such as Ackermann's, are also little-suited to investigate the notion of normal sets and to explicate them properly. Furthermore, Ackermann's system generates essentially only the cumulative hierarchy of ZF and hence can only yield results which we already know anyway.[21]

In NAM (and its modifications) arbitrarily many formulas can be taken as starting points for investigating the normality of sets. Some formulas will be found to be quite suitable. Usually all such formulas merely imply the predicate "Normal". But once such a predicate is found to be equivalent to the atom "Normal(x)" (which does not make the considered systems of NAM inconsistent), then the concept of normal set would be more or less correctly and uniquely characterized. (Of course, an entire hierarchy of characterizations with increasing strength is possible, and the most counterintuitive ones are excluded.) The same happens with NACT and the explication of the notion of pathological classes.

The approach suggested here is similar to the attempts made in the 30's to find an adequate explicatum or surrogate for the notion of computability. In NAM and NACT (and its modifications) one can experiment and study the effects of varying the normality conditions, until one has found the most suitable. With such normality conditions mathematicians will finally obtain a decision criterion for the generation of such sets in naïve set theory that they consider as most natural. I think this amounts to an important philosophical program.

## References for the Set-Theoretical Articles, K50-Set10

---

[21] See Felgner [2002] and K50-Set3 in this book.